\documentclass[12pt]{article}
\usepackage{amsmath,amssymb,theorem}
\usepackage{graphicx}
\usepackage{subfigure}
\usepackage{psfrag}
\usepackage{color}
\usepackage{epstopdf}
\usepackage[T1]{fontenc}

\newtheorem{thm}{Theorem}[section]
\newtheorem{conj}[thm]{Conjecture}
\newtheorem{cor}[thm]{Corollary}
\newtheorem{lem}[thm]{Lemma}
\theorembodyfont{\rmfamily}
\def\pf{\bigskip\noindent {\bf Proof.}~~}

\def\dfn#1{{\sl #1}}

\def\less{\backslash}

\def\pf{\bigskip\noindent {\emph{Proof.}}~~}

\def\mytextindent#1{\indent\llap{#1\enspace}\ignorespaces}
\def\fattextindent#1{\indent\indent\llap{#1\enspace}\ignorespaces}
\def\myitem{\par\hangindent\parindent\mytextindent}
\def\myitemitem{\par\hangindent\parindent\fattextindent}

\title{Star chromatic index    of subcubic multigraphs}

\author{Hui Lei\footnotemark[1] \, and Yongtang Shi\thanks{Partially supported by the National
Natural Science Foundation of China and Natural Science Foundation of Tianjin (No.
17JCQNJC00300).} \\
Center for Combinatorics and LPMC\\
Nankai University\\
 Tianjin 300071, China\\
~\\
Zi-Xia Song \thanks{Corresponding author.  Email: Zixia.Song@ucf.edu}\\
Department  of Mathematics\\
 University of Central Florida\\
   Orlando, FL 32816, USA
}
\begin{document}
\maketitle
\begin{abstract}
The \dfn{star chromatic index}  of a multigraph $G$, denoted
$\chi'_{s}(G)$,  is the minimum number of colors needed to properly
color the edges of $G$ such that no path or cycle of length four is
bi-colored.  A multigraph $G$ is \dfn{star $k$-edge-colorable} if
$\chi'_{s}(G)\le k$.   Dvo\v{r}\'ak, Mohar and \v{S}\'amal [Star
chromatic index, {\it J. Graph Theory} {\bf 72} (2013), 313--326]
proved  that every  subcubic multigraph is star $7$-edge-colorable.
They  conjectured in the same paper that every  subcubic multigraph
should be star $6$-edge-colorable. In this paper, we first prove that
 it is NP-complete to determine whether $\chi'_s(G)\le3$ for an arbitrary graph $G$.
 This answers a question of Mohar.  We then establish  some
structure results on subcubic multigraphs $G$ with $\delta(G)\le2$ such that $\chi'_s(G)>k$
but $\chi'_s(G-v)\le k$ for any $v\in V(G)$, where $k\in\{5,6\}$.
We finally apply the structure results,  along with a simple
discharging method,  to prove that every subcubic multigraph $G$  is
star $6$-edge-colorable if $mad(G)<5/2$,  and star $5$-edge-colorable
if $mad(G)<24/11$, respectively, where $mad(G)$ is
the maximum average degree of  a multigraph $G$. This partially confirms the conjecture of  Dvo\v{r}\'ak, Mohar and \v{S}\'amal. \\

\end{abstract}

\section{Introduction}

\baselineskip 18pt 

All multigraphs  in this paper are finite and loopless; and all graphs  are finite and without loops or multiple edges.
Given a multigraph $G$, let $c: E(G)\rightarrow [k]$ be a proper
edge-coloring of $G$, where $k\ge1$ is an integer and $[k]:=\{1,2, \dots, k\}$. We say that
$c$ is a  \dfn{star $k$-edge-coloring} of $G$ if no path or cycle of
length four in $G$ is bi-colored under the coloring $c$; and  $G$ is
\dfn{star $k$-edge-colorable} if $G$ admits a star
$k$-edge-coloring. The \dfn{star chromatic index}  of $G$, denoted
by $\chi'_{s}(G)$,  is the smallest integer $k$ such that $G$ is
star $k$-edge-colorable.  The  chromatic index of $G$ is denoted by
$\chi'(G)$. As pointed out in  \cite{DMS2013}, the definition of
star edge-coloring of a graph $G$ is equivalent to the star
vertex-coloring of its line graph $L(G)$.  Star edge-coloring of a
graph was initiated  by Liu and Deng \cite{DL2008}, motivated by the
vertex version (see \cite{ACKKR2004, BCMRW2009, CRW2013, KKT2009,
NM2003}).   Given a multigraph $G$, we use  $|G|$ to denote the number of vertices,  $e(G)$ the number of edges, $\delta(G)$ the minimum degree, and $\Delta(G)$ the maximum degree of $G$, respectively.  For any $v\in V(G)$, let $d_G(v)$  and $N_G(v)$ denote the degree   and neighborhood of $v$ in  $G$, respectively.  For any subsets  $A, B\subseteq V(G)$, let $N_G(A):=\bigcup_{a\in A} N_G(a)$, and    let  $A\less B := A-B$.   If  $B=\{b\}$, we simply write  $A\less b$ instead of $A\less B$.  We use $K_n$ and $P_n$ to denote the complete graph and the path on $n$ vertices, respectively. \medskip

It  is well-known~\cite{Vizing}  that the chromatic index of a graph with maximum degree $\Delta$ is   either $\Delta$ or $\Delta+1$. However,  it is NP-complete~\cite{Holyer} to determine whether the chromatic index of an
arbitrary graph with maximum degree $\Delta$ is $\Delta$ or $\Delta+1$.  The problem remains NP-complete even for cubic graphs.    A multigraph $G$ is \dfn{subcubic} if the maximum degree
of $G$ is at most three. Mohar (private communication with the second author) proposed that   it is  NP-complete to determine whether $\chi'_s(G)\le3$ for an arbitrary graph $G$.  We first  answer this question in the positive.\medskip

\begin{thm}\label{NP}
It is  NP-complete to determine whether $\chi'_s(G)\le3$ for an arbitrary graph $G$.
\end{thm}
\medskip

 We prove Theorem~\ref{NP} in Section~\ref{NPC}. Theorem~\ref{Kn}  below is a result of Dvo\v{r}\'ak, Mohar and \v{S}\'amal~\cite{DMS2013},   which gives  an   upper bound and a lower  bound for complete graphs.

\begin{thm}[\cite{DMS2013}]\label{Kn}  The star chromatic index of   the complete graph $K_n$ satisfies

$$2n(1+o(1))\leq \chi'_{s}(K_n)\leq n\, \frac{2^{2\sqrt{2}(1+o(1))\sqrt{\log n}}}{(\log n)^{1/4}}.$$
In particular, for every $\epsilon>0$, there exists a constant $c$ such that  $\chi'_{s}(K_n)\le cn^{1+\epsilon}$ for every integer $n\ge1$.
\end{thm}

The true order of magnitude of  $\chi'_{s}(K_n)$ is still unknown. From Theorem~\ref{Kn},  an upper bound in terms of the maximum degree for general  graphs is also derived in~\cite{DMS2013}, i.e.,
$\chi'_{s}(G)\leq \Delta\cdot 2^{O(1)\sqrt{\log \Delta}}$ for  any  graph $G$ with maximum degree $\Delta$.  In the same paper, Dvo\v{r}\'ak, Mohar and \v{S}\'amal~\cite{DMS2013} also considered the star chromatic index of subcubic multigraphs. To state their result, we need to introduce one notation. A graph $G$ \dfn{covers} a graph $H$ if there is a mapping $f: V(G)\rightarrow V(H)$ such that for any $uv\in E(G)$, $f(u)f(v)\in E(H)$, and for any $u\in V(G)$, $f$ is a bijection between $N_G(u)$ and $N_{H}(f(u))$.  They proved the following.

\begin{thm} [\cite{DMS2013}]\label{s=7} Let $G$ be a multigraph.\medskip

\myitemitem{(a)}  If $G$ is  subcubic,  then $\chi'_s(G)\le7$.\medskip

\myitemitem{(b)}  If $G$ is  cubic and has no multiple edges, then $\chi'_s(G)\ge4$ and the equality holds if and only if $G$ covers the graph of $3$-cube.

\end{thm}

As observed in~\cite{DMS2013},  $\chi'_s(K_{3,3})=6$ and the Heawood graph is star $6$-edge-colorable. No subcubic multigraphs with star chromatic index seven  are known. Dvo\v{r}\'ak, Mohar and \v{S}\'amal~\cite{DMS2013}   proposed the following conjecture.

\begin{conj}\label{cubic}
Let $G$ be a subcubic multigraph. Then $\chi'_s(G)\leq 6$.
\end{conj}

As far as we know, not much progress has been made yet  towards Conjecture~\ref{cubic}.  It was recently shown in~\cite{BLM2016} that every subcubic outerplanar graph is star $5$-edge-colorable.
A tight upper bound for trees was also obtained in  \cite{BLM2016}. We summarize   the main results in   \cite{BLM2016}  as follows.\medskip

\begin{thm}[\cite{BLM2016}]\label{outplanar} 
Let $G$ be an outerplanar graph. Then\medskip

\myitemitem{(a)} $\displaystyle\chi'_{s}(G)\leq\left\lfloor\frac{3\Delta(G)}2\right\rfloor$ if $G$ is a tree.  Moreover,  the bound is tight. \bigskip

\myitemitem {(b)} $\chi'_s(G)\le5$ if $\Delta(G)\le3$.\bigskip

\myitemitem {(c)}  $\displaystyle \chi'_s(G) \le\left\lfloor\frac{3\Delta(G)}{2}\right\rfloor+12$  if $\Delta(G)\ge4$. \medskip
\end{thm}

The \dfn{maximum average degree} of  a multigraph $G$, denoted  $mad(G)$, is defined as the maximum  of  $2
e(H)/|H|$ taken over all the subgraphs $  H$ of $G$.
 We want to point out here that there is an error in the proof of
Theorem 2.3 in a recent published paper by Pradeep and Vijayalakshmi
[Star chromatic index of subcubic graphs, {\it Electronic Notes in Discrete
Mathematics}  {\bf 53} (2016), 155--164]. Theorem 2.3 in~\cite{PV2016} claims  that if
$G$ is a subcubic  graph with $mad(G)<11/5$, then $\chi'_s(G)\le5$.
 The error  in the proof of
Theorem 2.3 arises from ambiguity in the
statement of Claim 3 in their paper. From its proof given in~\cite{PV2016} (on page 158),
Claim 3 should be stated as ``$H$ does not
 contain a path $uvw$,  where either all of $u, v,  w$ are $2$-vertices or all of $u,v,  w$ are light $3$-vertices".  This new statement of Claim 3
does not imply that ``a $2$-vertex must be adjacent to a heavy
$3$-vertex" in Case 2 of the proof of Theorem 2.3 (on page 162).  It seems nontrivial to fix this error in their proof. If  Claim 3 in their paper is  true, using the technique we developed in the proof of Theorem~\ref{thm*}(b), one can obtain  a stronger result that every subcubic multigraph with    $mad(G)<7/3$ is star $5$-edge-colorable. \medskip

In this paper, we  prove two main results, namely Theorem~\ref{NP} mentioned above  and Theorem~\ref{thm*} below.

\begin{thm}\label{thm*}
Let $G$ be a subcubic multigraph.\medskip

\myitemitem {(a)}   If $mad(G)<2$, then $\chi'_{s}(G)\leq 4$ and the bound is tight. \medskip

\myitemitem {(b)}  If $mad(G)<24/11$, then $\chi'_{s}(G)\leq 5$.\medskip

\myitemitem {(c)} If $mad(G)<5/2$, then $\chi'_{s}(G)\leq 6$.
\end{thm}
\medskip

The rest of this paper is organized as follows. We prove Theorem~\ref{NP} in Section~\ref{NPC}. Before we prove Theorem~\ref{thm*} in  Section~\ref{k=56}, we establish in Section~\ref{prop} some structure results on subcubic multigraphs $G$ with $\delta(G)\le2$ such that   $\chi'_s(G)>k$ and $\chi'_s(G-v)\le k$ for any $v\in V(G)$, where $k\in\{5,6\}$. We believe that our structure results can be used to solve Conjecture~\ref{cubic}. \medskip

\section{Proof of Theorem \ref{NP}}\label{NPC}
First let us denote by SEC the  problem stated in Theorem~\ref{NP}, and   we denote by  3EC the following well-known NP-complete problem of Holyer~\cite{Holyer}: \medskip

Given a  cubic graph $G$, is $G$ $3$-edge-colorable? \bigskip

\noindent {\bf Proof of  Theorem~\ref{NP}}:   Clearly, SEC is  in the class   NP.  We shall reduce 3EC  to SEC. \medskip

\begin{figure}[htbp]
\begin{center}
\scalebox{1.2}[1.2]{\includegraphics{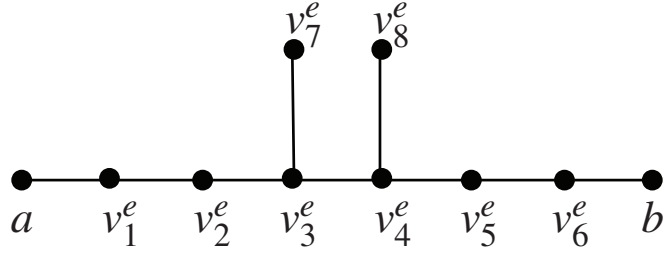}}
\caption{Graph $H_{ab}$.}\label{fig1}
\end{center}
\end{figure}

Let  $H$ be an instance of 3EC. We construct a   graph $G$ from $H$ by replacing   each edge $e=uw\in E(H)$  with a copy of  graph $H_{ab}$, identifying $u$ with $a$ and $w$ with $b$, where $H_{ab}$ is depicted in Figure \ref{fig1}.
The size of $G$ is clearly polynomial in the size of $H$, and $\Delta(G)=3$.

It suffices to show that   $\chi' (H)\le3$ if and only if $\chi'_{s}(G)\le3$.  Assume that $\chi'(H)\le3$.  Let $c:E(H)\rightarrow \{1,2,3\}$
be a proper $3$-edge-coloring of $H$. Let  $c^*$  be an edge coloring of $G$ obtained from $c$  as follows:
for each edge $e=uw\in E(H)$, let  $c^*(av^e_1)=c^*(v^e_3v^e_4)=c^*(v^e_6b)=c(uw)$,
$c^*(v^e_1v^e_2)=c^*(v^e_3v^e_7)=c^*(v^e_4v^e_8)=c^*(v^e_5v^e_6)=c(uw)+1$,
and $c^*(v^e_2v^e_3)=c^*(v^e_4v^e_5)=c(uw)+2$, where all colors here and henceforth are done modulo $3$.
Notice  that $c^*$ is a proper $3$-edge-coloring of $G$. Furthermore,  it can be easily checked  that
$G$ has no bi-colored path or  cycle of length four under the coloring $c^*$. Thus  $c^*$ is a star $3$-edge-coloring of $G$ and so  $\chi'_{s}(G)\le3$.\medskip

Conversely, assume that  $\chi'_{s}(G)\le3$. Let  $c^*:E(G)\rightarrow \{1,2,3\}$
be a star $3$-edge-coloring of $G$.  Let $c$ be an edge-coloring of $H$ obtained from $c^*$ by letting $c(e)=c^*(av^e_1)$ for any $e=uw\in E(H)$.  Clearly, $c$ is a proper $3$-edge-coloring of $H$ if  for any edge  $e=uw$ in $G$, $c^*(av^e_1)=c^*(v^e_6b)$. We prove this next.   Let $e=uw$ be an edge of $H$. We consider the following two cases. \medskip

\noindent {\bf Case 1:} $c^*(v^e_3v^e_7)=c^*(v^e_4v^e_8)$.\medskip

In this case, let  $c^*(v^e_3v^e_7)=c^*(v^e_4v^e_8)=\alpha$, where $\alpha\in\{1,2,3\}$. We may further assume that  $c^*(v^e_3v^e_4)=\beta$ and $c^*(v^e_2v^e_3)=c^*(v^e_4v^e_5)=\gamma$, where $\{\beta, \gamma\}=\{1,2,3\}\less \alpha$. This is possible because   $d_{G^*}(v^e_3)=d_{G^*}(v^e_4)=3$  and $c^*$ is a proper $3$-edge-coloring of $G^*$.  Since $c^*$ is a star edge-coloring of $G^*$, we see that $c^*(v^e_1v^e_2)=c^*(v^e_5v^e_6)=\alpha$ and so $c^*(av^e_1)=c^*(v^e_6b)=\beta$. \\

\noindent {\bf Case 2:} $c^*(v^e_3v^e_7)\ne c^*(v^e_4v^e_8)$.\medskip

In this case, let  $c^*(v^e_3v^e_7)=\alpha$, $c^*(v^e_4v^e_8)=\beta$, $c^*(v^e_3v^e_4)=\gamma$, where $\{\alpha, \beta, \gamma\}=\{1,2,3\}$. This is possible because $\alpha\ne \beta$ by assumption.
Since $c^*$ is a proper edge-coloring of $G^*$, we see that   $c^*(v^e_2v^e_3)=\beta$  and $c^*(v^e_4v^e_5)=\alpha$. One can easily check now that $c^*(v^e_1v^e_2)=\alpha$ and $c^*(v^e_5v^e_6)=\beta$,  and so $c^*(av^e_1)=c^*(v^e_6b)=\gamma$, because $c^*$ is a star edge-coloring of $G^*$. \medskip

In both cases we see that $c^*(av^e_1)=c^*(v^e_6b)$. Therefore  $c$ is  a proper $3$-edge-coloring of $H$ and so $\chi'(H)\le3$.  This completes the proof of Theorem~\ref{NP}.
\hfill\vrule height3pt width6pt depth2pt\\


\section{Properties of star $k$-critical subcubic multigraphs }\label{prop}

In this section,  we establish some structure results on subcubic multigraphs $G$ with $\delta(G)\le2$ such that   $\chi'_s(G)>k$ and  $\chi'_s(G-v)\le k$ for any $v\in V(G)$, where $k\in\{5,6\}$.  For simplicity, we say that  a  multigraph $G$ is  \dfn{star $k$-critical} if $\chi'_s(G)>k$ and  $\chi'_s(G-v)\le k$ for any $v\in V(G)$, where $k\in\{5,6\}$.
 Clearly, every star $k$-critical graph  must be connected.
\medskip

Throughout the remainder  of this section, let $G$ be a star $k$-critical subcubic multigraph with $\delta(G)\le2$, and let $N(v)$ and $d(v)$ denote the neighborhood and  degree of a vertex $v$ in $G$, respectively.  Since every multigraph with maximum degree two is star $4$-edge-colorable, we see that  $\Delta(G)=3$ and $|G|\ge3$.    Let $x\in V(G)$ with $d(x)\le 2$.  Let $H=G-x$ and let $c: E(H)\rightarrow [k]$ be a star $k$-edge-coloring of $H$, where $k\in\{5,6\}$.
For any $u\in V(H)$, let $c(u)$ denote the set of all colors such that each is used to color an edge incident with $u$ under the coloring $c$.
 For any $A\subseteq V(H)$, let $c(A):= \bigcup_{a\in A} c(a)$. By abusing the notation we use    $c(uv)$  to denote the set of all colors on the edges between $u$ and $v$ under the coloring $c$ if  $uv\in E(H)$  is a multiple edge. \bigskip

 \noindent {\bf Observation} \,\, If  $d(x)=2$, then $|N(x)|=2$.

\pf  Suppose that $|N(x)|=1$.  Let $N(x)=\{z\}$. Since $G$ is connected, we see that $d(z)=3$. Let $N(z)=\{x, z^*\}$. We obtain a star $k$-edge-coloring of $G$ by coloring the two edges between $x$ and $z$ by two distinct colors in $[k]\less c(z^*)$, a contradiction.   \hfill\vrule height3pt width6pt depth2pt\medskip

\begin{lem}\label{deg=1}
Assume that  $d(x)=1$. Let  $N(x)=\{y\}$. The following are true.\medskip

\myitem{(a)}  $c(N_H(y))=[k]$ and  $|N(y)|=3$.
\myitem{(b)} $N(y)$ is an independent set in $G$,  $d(y_1)=3$ and $d(y_2)\ge k-3$, where $N(y)=\{x, y_1,y_2\}$ with   $d(y_1)\ge d(y_2)$.
\myitem{(c)}  If $d(y_2)=k-3$, then for any $i\in\{1,2\}$ and any $v\in N_H(y_i)\less y$,  $c(yy_i)\in c(v)$, $|N(v)|\ge2$, $|N(y_1)|=3$, $|N(y_2)|=k-3$,  and $N[y_1]\cap N[y_2]=\{y\}$.
\myitem{(d)} If $d(y_2)=2$, then $k=5$ and $d(w_1)=3$, where $w_1$ is the other neighbor of $y_2$ in $G$.
\myitem{(e)} If  $k=6$ and  for some $i\in\{1,2\}$,  $N(y_i)\less y$  has a vertex $v$ with $d(v)=2$,  then $vv'\notin E(G)$, $N(v')$ is an independent set in $G$,  and $d(u)=3$ for any $u\in N(v)\cup N[v']$, where $N(y_i)=\{y, v, v'\}$.
\myitem{(f)} If  $k=6$ and  for some $i\in\{1,2\}$, each vertex of $N(y_i)\less y$ has degree three in $G$, then either $d(v)\ge2$ for any $v\in N(y^1_i)$ or $d(v)\ge2$ for any $v \in N(y^2_i)$, where $N(y_i)=\{y, y^1_i, y^2_i\}$.
\myitem{(g)} If $k=5$ and $d(y_2)=3$, then either $d(v)\ge2$ for any $v\in N(y_1)$ or $d(v)\ge2$ for any $v \in N(y_2)$.
\end{lem}
\pf To prove Lemma~\ref{deg=1}(a), suppose that $c(N_H(y))\ne [k]$. Then coloring  the edge $xy$ by a color in $[k]\less c(N_H(y))$, we obtain a star $k$-edge-coloring of $G$, a contradiction. Thus  $c(N_H(y))=[k]$ and so  $|N(y)|=3$. \medskip

Next let $N(y)=\{x, y_1,y_2\}$ with   $d(y_1)\ge d(y_2)$ by Lemma~\ref{deg=1}(a). Suppose that  $y_1y_2\in E(G)$.  Then  $|N(y_1)|=|N(y_2)|=3$ and  all the edges incident to $y_1$ or $y_2$ are colored with distinct colors because $c(N_H(y))= [k]$.  Now coloring the edge $xy$ by color $c(y_1y_2)$ yields a star $k$-edge-coloring of $G$, a contradiction. Thus $y_1y_2\notin E(G)$.  Since $|c(N_H(y))|= k\ge5$, we see that   $d(y_1)=3$ and $d(y_2)\ge k-3\ge2$. This proves Lemma~\ref{deg=1}(b).\medskip

 To prove Lemma~\ref{deg=1}(c), since  $d(y_2)= k-3$, we see that  all the edges incident to $y_1$ or $y_2$ are colored with distinct colors because $c(N_H(y))= [k]$.   Suppose that for some $i\in\{1,2\}$, there exists a vertex $v\in N_H(y_i)\less y$ such that $c(yy_i)\notin c(v)$. Then we obtain a star $k$-edge-coloring of $G$ by coloring the edge $xy$ by a color in $c(y_iv)$, a contradiction. Thus for any $i\in\{1,2\}$ and any $v\in N_H(y_i)\less y$,  $c(yy_i)\in c(v)$. Hence $|N(v)|\ge2$.  Since  $\Delta(G)=3$, we see that $N[y_1]\cap N[y_2]=\{y\}$.  We next show that
  $|N(y_1)|=3$ and $|N(y_2)|=k-3$. \medskip

  Suppose that $|N(y_1)|<3$. Then $|N(y_1)|=2$ because  $d(y_1)=3$. Let $N(y_1)=\{y, y_1^*\}$. Then $d(y_1^*)=3$ by Observation.  Let $u$ be the other neighbor of $y_1^*$ in $G$.
 Then $c(y_1^*u)=c(yy_1)$. Let $\alpha, \beta$ be  two distinct colors on the parallel edges  $y_1y_1^*$. Then $\alpha, \beta\notin c(u)$ because $c$ is a star $k$-edge-coloring of $H$.
 Let $e^*$ be the edge between $y_1$ and $y_1^*$ with color $\alpha$. If $c(y_2)\less c(u)\ne\emptyset$, then we obtain a star $k$-edge-coloring of $G$ by recoloring the edge $e^*$ by a color in  $c(y_2)\less c(u)$, $yy_1$ by color $\alpha$, and coloring the edge $xy$ by color $c(yy_1)$, a contradiction. Thus $c(y_2)\subset  c(u)$ and so $k=5$. Clearly, $|N(y_2)|=k-3=2$. Let $y_2^*$ be the other neighbor of $y_2$ in $G$. Then $c(yy_2) \in c(y^*_2)$.  We obtain a star $5$-edge-coloring of $G$ by recoloring the edge $yy_2$ by a color in  $\{\alpha, \beta\}\less c(y_2^*)$, $yy_1$ by color $c(yy_2)$, and coloring the edge $xy$ by color $c(yy_1)$, a contradiction. Thus  $|N(y_1)|=3$. By symmetry, $|N(y_2)|=3$ if $k=6$. Clearly, $|N(y_2)|=k-3=2$ if $k=5$.  This proves Lemma~\ref{deg=1}(c). \medskip

It remains to prove Lemma~\ref{deg=1}(d), (e), (f) and (g).  Notice that if $d(y_2)=2$, then by Lemma~\ref{deg=1}(b),  $d(y_2)=2\ge k-3$. Thus $k=5$ and  so  $d(y_2)=k-3$. For each proof of  Lemma~\ref{deg=1}(d), (e), and (f), let  $d(y_2)=k-3$.  By Lemma~\ref{deg=1}(c),  we may assume that $N(y_1)=\{y, z_1, z_2\}$ and  $N(y_2)=\{y, w_1\}$ if $k=5$ (and $N(y_2)=\{y, w_1, w_2\}$ if $k=6$), where $N[y_1]\cap N[y_2]=\{y\}$.  By  Lemma~\ref{deg=1}(a), we may further assume that $c(yy_1)=1$, $c(yy_2)=2$, $c(y_1z_1)=3$,  $c(y_1z_2)=4$,  $c(y_2w_1)=5$  (and $c(y_2w_2)=6$ when $k=6$).
By Lemma~\ref{deg=1}(c), $1\in c(z_1)\cap c(z_2)$ and $2\in c(w_1)$ if $k=5$  (and $2\in c(w_1)\cap c(w_2)$ if $k=6$). \medskip

We next prove Lemma~\ref{deg=1}(d). Clearly, $k=5$.  Suppose that $d(w_1)\le2$.  By Lemma~\ref{deg=1}(c), $d(w_1)=2$  and $c(w_1)=\{2,5\}$.    Let $w^*\in N(w_1)$ with $c(w_1w^*)=2$. Notice that $w^*$ is not necessarily  different from $z_1$ or $z_2$.    Since $c$ is a star edge-coloring of $H$, $5\notin c(w^*)$. If $3\notin c(w^*)$, then we obtain a star $5$-edge-coloring of $G$ by coloring the edge $xy$ by color $5$ and  recoloring the edge $y_2w_1$  by color $3$, a contradiction. Thus $3\in c(w^*)$, and similarly, $4\in c(w^*)$. Hence  $w^*\notin \{z_1, z_2\}$ because $\Delta(G)=3$. We obtain a star $5$-edge-coloring of $G$ by coloring the edge $xy$ by color $2$,   recoloring the edge $yy_2$ by color $5$, and $y_2w_1$ by  color $1$, a contradiction. \medskip

To prove Lemma~\ref{deg=1}(e), since $k=6$, we see that $c(y_1)=\{1,3,4\}$ and  $c(y_2)=\{2,5,6\}$.  By symmetry, we may assume that $i=1$,   $v=z_1$,  and $v'=z_2$.   Then $1\in c(z_1)\cap c(z_2)$.    Suppose that $z_1z_2\in E(G)$. Then $c(z_1z_2)=1$.  Now recoloring the edge  $y_1z_1$ by a color in $\{5,6\}\less c(z_2)$ and coloring  the edge $xy$ by color $3$, we obtain a star $6$-edge-coloring of $G$, a contradiction. Thus $z_1z_2\notin E(G)$.  Let  $N(z_1)=\{y_1, z_1^*\}$. Then $c(z_1z_1^*)=1$. Since $z_1z_2\notin E(G)$, we see that $z_1^*\ne z_2$.  If $\{2,5,6\}\less (c(z_1^*)\cup c(z_2))\ne\emptyset$, then recoloring  the edge $y_1z_1$ by a color in $\{2,5,6\}\less (c(z_1^*)\cup c(z_2))$,  $y_1y$ by color $3$,  and coloring $xy$ by color $1$ yields a star $6$-edge-coloring of $G$, a contradiction. Thus $\{2,5,6\}\subset c(z_1^*)\cup c(z_2)$ and so $d(z_1^*)=d(z_2)=3$. Clearly, $z_1^*z_2\notin E(G)$ because $\Delta(G)=3$. Let $c(z_2)=\{1,4,\alpha\}$, where $\alpha\in\{2,5,6\}$.
Suppose that $c(N(z_2)\less y_1)\ne [6]$. Then we obtain a star $6$-edge coloring of $G$ by recoloring the edge $y_1z_2$ by a color, say $\beta$,  in $[6]\less c(N(z_2)\less y_1)$, $y_1z_1$ by color $\alpha$, $yy_1$ by  color  $3$ if $\beta\ne3$ or color $4$ if $\beta=3$,
and finally coloring $xy$ by color $1$, a contradiction.
Thus $c(N(z_2)\less y_1)=[6]$ and so $|N(z_2)|=3$. Let $N(z_2)=\{y_1, z_2^1, z_{2}^2\}$. Clearly, $z_2^1z_2^2\notin E(G)$ because $c(z_2^1)\cup c(z_2^2)=[6]$, and $d(z_2^1)=d(z_2^2)=3$, as desired.  This proves Lemma~\ref{deg=1}(e).  \medskip

To prove Lemma~\ref{deg=1}(f), we may assume that $i=1$, $z_1=y^1_1$ and $z_2=y^2_1$.  Suppose that there exist  vertices  $z_1^1\in N(z_1)$ and  $z_2^1\in  N(z_2)$  such that $d(z_1^1)=d(z_2^1)=1$. Then $|N(z_1)|=|N(z_2)|=3$ by Lemma~\ref{deg=1}(a). Let $N(z_1)=\{y_1, z_1^1, z_{1}^2\}$ and $N(z_2)=\{y_1, z_2^1, z_{2}^2\}$. Then $d(z_1^2)=d(z_2^2)=3$ by Lemma~\ref{deg=1}(a).  Assume first that $c(z_1z_1^1)=1$.  If $\{2,5,6\}\less c(z_1^2)\ne\emptyset$, then recoloring the edge $z_1z_1^1$ by a color in $\{2,5,6\}\less c(z_1^2)$,  we obtain a star $6$-edge-coloring of $H$ with $1\notin c(z_1)$, contrary to Lemma~\ref{deg=1}(c). Thus $c(z_1^2)=\{2,5,6\}$. Clearly,  $3\in c(z_2)$, otherwise recoloring the edge $z_1z_1^1$ by  color $4$ yields  a star $6$-edge coloring of $H$ with $1\notin c(z_1)$, contrary to Lemma~\ref{deg=1}(c). Thus $c(z_2)=\{1,3,4\}$.  Then
$\{2,5,6\}\less c(z_2^2)\ne\emptyset$. Now recoloring the edge $z_2z_2^1$ by a color in $\{2,5,6\}\less c(z_2^2)$,  we obtain a star $6$-edge-coloring of $H$ with $c(z_2)\ne \{1,3,4\}$,
a contradiction. Thus  $c(z_1z_1^1)\ne1$. By symmetry, $c(z_2z_2^1)\ne1$. By Lemma~\ref{deg=1}(c), $c(z_1z_1^2)=c(z_2z_2^2)=1$. Clearly, $3\notin c(z_1^2)$
because $H$ has no bi-colored path of length four. If $\{5,6\}\less c(z_1^2)\ne \emptyset$, then we obtain a star $6$-edge-coloring of $G$ by recoloring the edge $z_1z_1^1$ by color $3$, $z_1y_1$ by a color in $\{5,6\}\less c(z_1^2)$,  and coloring $xy$ by color $3$,  a contradiction.  Thus $c(z_1^2)=\{1,5,6\}$. Similarly, $c(z_2^2)=\{1,5,6\}$. Now recoloring the edges $z_1z_1^1$ by color $3$, $z_1y_1$ by color $4$, $z_2^1z_2$ by color $4$, $y_1z_2$ by color $2$, $yy_1$ by color $3$, and finally coloring the edge $xy$ by color $1$, we obtain a star $6$-edge-coloring of $G$, a contradiction. This proves Lemma~\ref{deg=1}(f). 
\medskip

It remains to prove Lemma~\ref{deg=1}(g). Suppose that there exist  vertices  $y_1^1\in N(y_1)$ and  $y_2^1\in  N(y_2)$  such that $d(y_1^1)=d(y_2^1)=1$. Then $|N(y_1)|=|N(y_2)|=3$ by Lemma~\ref{deg=1}(a).  Let $N(y_1)=\{y, y_1^1, y_1^2\}$ and $N(y_2)=\{y, y_2^1, y_2^2\}$. By Lemma~\ref{deg=1} (a),  $c(y_1)\cup c(y_2)=[5]$. If $c(y_1y_1^1)=c(y_2y_2^1)$, then  we obtain a star $5$-edge-coloring of $G$ by coloring the edge $xy$ by color $c(y_1y_1^1)$ because $c(y_1)\cup c(y_2)=[5]$, a contradiction.  Thus $c(y_1y_1^1)\ne c(y_2y_2^1)$. Since $c(y_1)\cup c(y_2)=[5]$, we see that  either $c(y_1y_1^1)\notin c(y_2)$ or  $c(y_2y_2^1)\notin c(y_1)$.   We may assume that $c(y_1y_1^1)\notin c(y_2)$. But then  coloring the edge $xy$ by color $c(y_1y_1^1)$ yields a star $5$-edge-coloring of $G$, a contradiction. \medskip


This completes the proof of Lemma~\ref{deg=1}. \hfill\vrule height3pt width6pt depth2pt\medskip


 \begin{lem}\label{deg=2}
Assume that $d(x)=2$. Let  $N(x)=\{z, w\}$ with $|N(z)|\le |N(w)|$.\medskip

\myitem{(a)}  If $zw\in E(G)$, then  $k=5$, $|N(z)|=|N(w)|=3$ and $d(v)\ge2$ for any $v\in N(z)\cup N(w)$.
\myitem{(b)} If  $zw\notin E(G)$, then  $|N(w)|=3$ or  $k=5$,  $ |N(w)|=|N(z)|=2$,  and $d(w)=d(z)=3$.
\myitem{(c)}  If $d(z)=2$ and  $z^*w\in E(G)$,  then $k=5$,  $|N(z^*)|=|N(w)|=3$, and $d(u)=3$ for any $u\in (N[w] \cup N[z^*])\less \{x,z\}$, where $z^*$ is the other neighbor of $z$ in $G$.
\myitem{(d)} If  $k=6$ and $d(z)=2$, then  $z^*w\notin E(G)$,  $|N(z^*)|=|N(w)|=3$, and for any $v\in (N(w)\cup N(z^*))\less \{x,z\}$,   $d(v)=3$ and $d(u)\ge2$ for any $u\in N(v)$, where $N(z)=\{x, z^*\}$.
\myitem{(e)} If $k=5$ and $d(z)=2$, then  $|N(z^*)|=|N(w)|=3$, and $|N(v)|\geq2$ for any $v\in N(w)\cup N(z^*)$, where $N(z)=\{x, z^*\}$.
\end{lem}

\pf   Assume that $zw\in E(G)$. Since $G$ is connected, we see that $|N(w)|=3$. Let $N(w)=\{x, z, w^*\}$. We first show that  $|N(z)|=3$ and $N(z)\cap N(w)=\{x\}$. Suppose that
 $|N(z)|=2$ or $|N(z)|=3$ and $zw^*\in E(G)$. Then $|c(w)\cup c(w^*)|\le4$ when $c(zw)\notin c(w^*)$ and $|c(w) \cup  c(w^*)|\le 3$ when  $c(zw)\in c(w^*)$. We obtain a star $k$-edge-coloring of $G$ by coloring  the edge $xw$ by a color, say $\alpha$,  in $[k]\less (c(w)\cup c(w^*))$ and  then coloring $xz$ by color $c(ww^*)$ if $c(zw)\notin c(w^*)$ or a color in $[k]\less (c(w)\cup c(w^*)\cup\{\alpha\})$ if $c(zw)\in c(w^*)$, a contradiction.  Thus  $|N(z)|=3$ and  $z^*\ne w^*$, where  $N(z)=\{x, w, z^*\}$.  We next show that $k=5$. Suppose that $k=6$. Then $c(zw)\in c(z^*)$, otherwise  we obtain a star $6$-edge-coloring of $G$ by coloring  the edge $xz$ by a color, say $\alpha$,  in $[6]\less (c(z^*)\cup c(w))$ and  $xw$ by a color in $[6]\less(c(w)\cup c(w^*)\cup\{\alpha\})$.    We then obtain a star $6$-edge-coloring of $G$ by coloring the edge $xw$ by a color, say $\beta$,   in $[6]\less (c(w^*)\cup c(z))$ and $xz$ by a color in $[6]\less (c(z^*)\cup \{c(ww^*), \beta\})$,  a contradiction. Hence $k=5$. By Lemma~\ref{deg=1}(b), we see that $d(v)\ge2$ for any $v\in N(z)\cup N(w)$. This proves Lemma~\ref{deg=2}(a). \medskip

To prove Lemma~\ref{deg=2}(b),  by Lemma~\ref{deg=1}(a), $|N(w)|\ge |N(z)|\ge 2$.  We are done if $|N(w)|=3$. So we may assume that  $|N(w)|=2$. Then $|N(z)|=2$ because $|N(z)|\le |N(w)|$. Let $z^*$ and $w^*$ be the other neighbor of $z$ and $w$, respectively. If $ww^*$ or $zz^*$ is not a multiple edge (say the former) or $k=6$, then we obtain a star $k$-edge-coloring of $G$ by coloring the edge $xz$ by a color, say $\alpha$, in $[k]\less (c(z^*)\cup\{c(ww^*)\})$ and $xw$ by a color in $[k]\less (c(w^*)\cup\{\alpha\})$, a contradiction. Thus $|c(ww^*)|=|c(zz^*)|=2$ and $k=5$.  We see that $d(w)=d(z)=3$. \medskip

 We next prove Lemma~\ref{deg=2}(c).  Since $d(z)=2$,  we see that $zw\notin E(G)$ by Lemma~\ref{deg=2}(a). Then  $|N(w)|=3$ by Lemma~\ref{deg=2}(b).  Since $xz^*\notin E(G)$, by Lemma~\ref{deg=2}(b) again,   $|N(z^*)|=3$.  Let $N(w)=\{x, z^*, w^*\}$ and $N(z^*)=\{z, z^*_1, w\}$.
We first show that $k=5$. Suppose that $k=6$.
 Then $c$ can be extended to be a star $6$-edge-coloring of $G$ by coloring the edge $xw$ by color $c(zz^*)$ if  $c(zz^*)\notin  c(w^*)$ or a color $\alpha$ in $[6]\less (c(z^*)\cup c(w^*))$ if $c(zz^*)\in  c(w^*)$, and then coloring the edge $xz$ by a color in $[6]\less (c(z^*)\cup \{\alpha, c(ww^*)\})$,  a contradiction.  Thus $k=5$. We next show that $d(w^*)=3$. Suppose that $d(w^*)\le2$.  If $c(zz^*)\notin c(w^*)$, then we obtain a star $5$-edge-coloring of $G$ by coloring the edge $xw$ by color $c(zz^*)$ and $xz$ by a color in $[5]\less (c(z^*)\cup c(w))$, a contradiction. Thus $c(zz^*)\in c(w^*)$ and so $|c(w^*)\cup c(z^*)|\le4$.  We obtain a star $5$-edge coloring of $G$ by coloring the edge $xw$ by a color, say $\alpha$,   in $[5]\less (c(z^*)\cup c(w^*))$ and $xz$ by a color in $[5]\less (c(z^*)\cup \{\alpha\})$, a contradiction.  By symmetry, $d(z^*_1)=3$. This proves Lemma~\ref{deg=2}(c).\medskip

To prove Lemma~\ref{deg=2}(d), since $k=6$, by Lemma~\ref{deg=2}(a,c), we see that  $wz, wz^*\notin E(G)$.  By  Lemma~\ref{deg=2}(b),    $|N(z^*)|=|N(w)|=3$.  Let $N(w)=\{x, w_1, w_2\}$ and $N(z^*)=\{z, z^*_1, z^*_2\}$. Then $c(zz^*)\in c(w_1)\cup c(w_2)$, otherwise we obtain a star $6$-edge-coloring of $G$ by coloring the edge $xw$ by color $c(zz^*)$ and $xz$ by a color in $[6]\less (c(z^*)\cup c(w))$, a contradiction.   We next show that $c(w_1)\cup c(w_2)=[6]$. Suppose that $c(w_1)\cup c(w_2)\ne [6]$.  Now coloring the edge $xw$ by a color, say $\alpha$,  in $[6]\less (c(w_1)\cup c(w_2))$,  and then coloring $xz$ by color $c(ww_1)$ if $c(ww_1)\notin c(z^*)$ or  a color  in $[6]\less (c(z^*)\cup c(w)\cup\{\alpha\})$ if $c(ww_1)\in c(z^*)$, we obtain  a star $6$-edge-coloring of $G$, a contradiction. Thus  $c(w_1)\cup c(w_2)=[6]$ and so $d(w_1)=d(w_2)=3$. By symmetry,  $d(z^*_1)=d(z^*_2)=3$.  Finally, for any $u\in N(z_1^*)\cup N(z_2^*)\cup N(w_1)\cup N(w_2)$, by Lemma~\ref{deg=1}(e),  $d(u)\ge2$. This proves Lemma~\ref{deg=2}(d).\medskip

It remains to prove Lemma~\ref{deg=2}(e). By Lemma~\ref{deg=2}(a), $zw\notin E(G)$. We may  assume that $wz^*\notin E(G)$ by Lemma~\ref{deg=2}(c). By Lemma~\ref{deg=2}(b),  $|N(z^*)|=|N(w)|=3$.  Suppose that  there exists a vertex $v\in N(w)\cup N(z^*)$ with $d(v)=1$. We may assume that $v\in N(w)$. Let $N(w)=\{x, v, w^*\}$.   Then  $c(zz^*)\in c(v)\cup c(w^*)$, otherwise we  recolor the edge $wv$ by color $c(zz^*)$. Now we  obtain a star $5$-edge-coloring of $G$ by coloring the edge $xw$ by  a color, say $\alpha$,  in $[5]\less (c(v)\cup c(w^*))$ and $xz$ by a color in $[5]\less (c(z^*)\cup \{\alpha\})$, a contradiction. Thus $|N(v)|, |N(w^*)|\ge2$. By symmetry, $|N(u)|\ge2$ for any $u\in N(z^*)$. \medskip

 This completes the proof of Lemma~\ref{deg=2}.  \hfill\vrule height3pt width6pt depth2pt\medskip

\begin{cor}\label{deg=3}
Let $G$ be a  subcubic multigraph that is star $6$-critical. Let $v\in V(G)$ be a vertex  with $N(v)=\{v_1, v_2, v_3\}$ and $d(v_1)\ge d(v_2)\ge d(v_3)=2$. 
The following are true. \medskip

\myitem {(a)}  $v_1, v_2, v_3$ are pairwise distinct, $d(v_1)=3$ and $v_3^*\notin\{v_1, v_2\}$, where $N(v_3)=\{v, v_3^*\}$.
\myitem {(b)} If   $d(v_3^*)=2$, then $d(v_2)=3$,  and  $d(u)\ge2$ for each $u\in N(v_1)\cup N(v_2)$,
\myitem {(c)}  If  $d(v_2)=2$, then every vertex of $N(v_2)\cup N(v_3)$ has degree three in $G$,  and for any $u\in N(v_1)\less v$,  $d(u)\ge2$.
\myitem {(d)}  If  $d(v_2)=3$ and there exists a vertex $v_i^*\in N(v_i)$ with $d(v_i^*)=1$ for some $i\in\{1,2\}$, then for any $u\in N(v_{3-i})$, $d(u)=3$.
\end{cor}

\pf   To prove Corollary~\ref{deg=3}(a), we first show that $v_1, v_2, v_3$ are pairwise distinct. By Observation,  $v_3$ is distinct from  $v_1, v_2$. Suppose that $v_1=v_2$. Let $v_1^*$  be the other neighbor of $v_1$, where $v_1^*,  v_3^*$ are not necessarily distinct. Let $c: E(G\less v)\rightarrow [6]$ be a star edge-coloring of $G\less v$. We obtain a star $6$-edge-coloring of $G$  by coloring the edges $vv_1$ by two distinct colors, say $\alpha, \beta$, in $[6]\less (c(v_1^*)\cup\{c(v_3v_3^*)\})$ and $vv_3$ by a color in $[6]\less(c(v_3^*)\cup\{\alpha, \beta\})$, a contradiction. This proves that $v_1, v_2, v_3$ are pairwise distinct. By Lemma~\ref{deg=2}(a), $v_3^*\notin\{v_1, v_2\}$. We next show that $d(v_1)=3$. Suppose that  $d(v_1)\ne3$. Then  $d(v_i)=2$ for all $i\in [3]$. By Lemma~\ref{deg=2}(a), we see that $v_1v_2\notin E(G)$. Let $v_i^*$ be the other neighbor of $v_i$ for all $i\in\{1,2\}$, where $v_1^*, v_2^*, v_3^*$ are not necessarily distinct. Let $c: E(G\less v)\rightarrow [6]$ be a star edge-coloring of $G\less v$. We obtain a star $6$-edge-coloring of $G$  by coloring the edge $vv_1$ by a color, say $\alpha$,  in $[6]\less (c(v_1^*)\cup\{c(v_2v_2^*), c(v_3v_3^*)\})$, $vv_2$ by a color, say $\beta$,  in $[6]\less (c(v_2^*)\cup\{\alpha, c(v_3v_3^*)\})$, and $vv_3$ by a color in $[6]\less (c(v_3^*)\cup\{\alpha, \beta\})$,  a contradiction. Thus $d(v_1)=3$.  This proves Corollary~\ref{deg=3}(a). \medskip

By Lemma~\ref{deg=2}(d),  Corollary~\ref{deg=3}(b) is true.  By Lemma~\ref{deg=2}(d) and Lemma~\ref{deg=1}(e),  Corollary~\ref{deg=3}(c) is true.  Finally, by Lemma~\ref{deg=1}(b,e) applied to $v_i^*$,   Corollary~\ref{deg=3}(d) is true. \hfill\vrule height3pt width6pt depth2pt\medskip

\section{Proof of Theorem \ref{thm*}}\label{k=56}

We are now ready to prove Theorem~\ref{thm*}.\medskip

\begin{figure}[htbp]
\begin{center}
\includegraphics[scale=0.8]{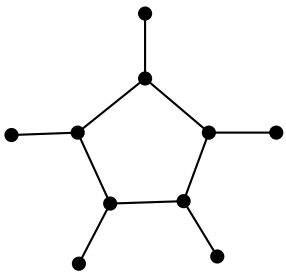}
\caption{A graph $G$ with $mad(G)=2$ and $\chi'_s(G)=5$.}\label{fig2}
\end{center}
\end{figure}

To prove Theorem~\ref{thm*}(a),  let $G$ be a subcubic multigraph with $mad(G)<2$. Then $G$ must be a simple graph. Notice that a simple graph $G$ has $mad(G)<2$
if and only if $G$ is a forest. Now applying  Theorem~\ref{outplanar}(a) to
every component of $G$, we see that $\chi'_s(G)\leq4$. This bound
is sharp in the sense that there exist  graphs $G$  with $mad\,
(G)=2$ and  $\chi'_s(G)>4$. One such
example from \cite{BLM2016} is depicted in Figure~\ref{fig2}.\medskip

We next proceed the proof of Theorem~\ref{thm*}(c) by contradiction.  Suppose the  assertion is false. Let $G$ be a subcubic multigraph with $mad(G)<5/2$ and  $\chi'_s(G)>6$. Among all counterexamples we choose $G$ so that
$|G|$ is minimum. By the choice of $G$, $G$ is connected,  star $6$-critical,  and  $mad(G)<5/2$. For all $i\in[3]$, let $A_i=\{v\in V(G): \, d_G(v)=i\}$ and  $n_i =|A_i|$ for all $i\in[3]$.
By Lemma~\ref{deg=1}(a), $A_1$  is an  independent set in $G$ and $N_G(A_1)\subseteq A_3$.  Let $G^*=G\less A_1$.  Then $mad(G^*)<5/2$.  We see that  $2e(G^*)=2e(G)-2n_1=2n_2+3n_3-n_1<5(n_2+n_3)/2$ and so  $n_3<n_2+2n_1$. Thus $A_1\cup A_2\ne\emptyset$. 
By Lemma~\ref{deg=1}(b), $\delta(G^*)\ge2$.   We say that a vertex $v\in V(G^*)$ with $d_{G^*}(v)=2$ is \dfn{good} if  $d_G(v)=3$; \dfn{bad} if $d_G(v)=2$ and $v$ is adjacent to another vertex of degree two in $G$; and \dfn{fair} if $d_G(v)=2$  and $v$ is not bad. We shall apply the discharging method  below to obtain a contradiction. \medskip

For each vertex $v\in V(G^*)$,   let  $\omega(v):= d_{G^*}(v)-\frac5 2$ be the initial charge of $v$.   Then $ \sum_{v\in V(G^*)} \omega(v) =2e(G^*)-\frac52|G^*|=(2n_2+3n_3-n_1)-\frac52(n_2+n_3)={(n_3-n_2-2n_1)}/2<0$,  because $n_3<n_2+2n_1$. Notice that for each $v\in V(G^*)$, $\omega(v)=2-5/2=-1/2$ if $d_{G^*}(v)=2$, and $\omega(v)=3-5/2=1/2$ if $d_{G^*}(v)=3$. Let $x\in V(G^*)$ be a vertex with $d_{G^*}(x)=3$ such that $x$  is  adjacent to  exactly $t\ge1$ vertices of degree two in $G^*$. We claim that $t\le2$ and $t=1$ when  $N_{G^*}(x)$ has a bad vertex. Clearly, $N_{G^*}(x)$ has at most two good vertices by Lemma~\ref{deg=1}(f), and at most two fair vertices by Corollary~\ref{deg=3}(a). By  Corollary~\ref{deg=3}(b), $N_{G^*}(x)$ has at most one bad vertex, and if such a bad vertex exists, then  $N_{G^*}(x)$ has no good or fair vertex.  Finally,  $N_{G^*}(x)$ has at most one fair vertex and one good vertex simultaneously by Corollary~\ref{deg=3}(c,d).  Thus $t\le2$ and $t=1$ when $N_{G^*}(x)$ has a bad vertex, as claimed. We will redistribute the charges of $x$  according to the following discharging rule: \medskip

\noindent{\bf  (R)}: For each $x\in V(G^*)$ with $d_{G^*}(x)=3$ and  exactly $t\ge1$ neighbors of degree two in $G^*$,   $x$ sends $\frac1{2t} \ge \frac14$ charges to each of its neighbors of degree two in $G^*$. \medskip

Let $\omega^*$ be the new charge of $G^*$ after applying the above discharging rule. We see that for any $v\in V(G^*)$ with $d_{G^*}(v)=3$, $\omega^*(v)\ge0$. We next show that for any $v\in V(G^*)$ with $d_{G^*}(v)=2$,
$\omega^*(v)\ge0$.
Let $v\in V(G^*)$ be a vertex with $d_{G^*}(v)=2$. By Observation and Lemma~\ref{deg=1}(a),   $|N_{G^*}(v)|=2$. Let $N_{G^*}(v)=\{u,w\}$. If $v$ is good, then $d_{G^*}(u)=d_{G^*}(w)=3$ by Lemma~\ref{deg=1}(c). Thus $\omega^*(v)\ge \omega(v)+1/4+1/4=0$.
Next, suppose that $v$ is bad. We may assume that $u$ is bad. By Lemma~\ref{deg=2}(d), $d_{G^*}(w)=3$.  By the above claim, $v$ is the only (bad) vertex of degree two of $N_{G^*}(w)$.
By the discharging rule (R),  $\omega^*(v)\ge \omega(v)+1/2=0$.  Finally, suppose that $v$ is a fair vertex.  Then $d_{G}(u)=d_{G}(w)=3$.  By Lemma~\ref{deg=1}(b),  neither $u$ nor $w$ is good. Thus $d_{G^*}(u)=d_{G^*}(w)=3$. By the discharging rule (R),
 $\omega^*(v)\ge \omega(v)+1/4+1/4=0$.  This proves that $\omega^*(v)\ge0$ for any $v\in V(G^*)$ with $d_{G^*}(v)=2$.  Thus  $\sum_{v\in V(G^*)} \omega^*(v) \ge0$, contrary to the fact that $ \sum_{v\in V(G^*)} \omega^*(v)=\sum_{v\in V(G^*)} \omega(v) <0$. \medskip

This completes the proof of Theorem~\ref{thm*}(c). \hfill\vrule height3pt width6pt depth2pt\medskip

The proof of Theorem~\ref{thm*}(b) is similar to the proof  Theorem~\ref{thm*}(c). For its completeness, we include its proof here because  the discharging part is different and more involved.   Suppose the  assertion is false. Let $G$ be a subcubic multigraph with $mad(G)<24/11$ and $G$ is not star $5$-edge-colorable. Among all counterexamples we choose $G$ so that
$|G|$ is minimum. By the choice of $G$, $G$ is connected and  star $5$-critical. Clearly,  $mad(G)<24/11$.  For all $i\in[3]$, let $A_i=\{v\in V(G): \, d_G(v)=i\}$ and let  $n_i=|A_i|$ for all $i\in[3]$.
 By Lemma~\ref{deg=1}(a), $A_1$  is an  independent set in $G$ and $N_G(A_1)\subseteq A_3$.  Let $G^*=G\less A_1$.   Then $mad(G^*)<24/11$.  We see that  $2e(G^*)=2e(G)-2n_1=2n_2+3n_3-n_1<24(n_2+n_3)/11$ and so  $9n_3<2n_2+11n_1$.  Thus $A_1\cup A_2\ne\emptyset$. 
 By Lemma~\ref{deg=1}(b), $\delta(G^*)\ge2$.    We say that a vertex $v\in V(G^*)$ with $d_{G^*}(v)=2$ is \dfn{good} if  $d_G(v)=3$;  \dfn{bad} if $d_G(v)=2$ and $v$ is adjacent to another vertex of degree two in $G$; and \dfn{fair} if $d_G(v)=2$  and $v$ is not bad.  Let $B:=\{v\in V(G^*): \, d_{G^*}(v)=2\}$.  We next claim  that every component of $G^*[B]$ is isomorphic to $K_1$, $K_2$ or $P_3$.  \medskip

 Suppose not. Let $P$ be a longest path in  $G^*[B]$ with vertices $x_1, x_2, \dots, x_p$ in order, where $p\ge4$ or $p=3$ and $x_1x_3\in E(G)$ or $p=2$ and $x_1x_2$ is a multiple edge.   Clearly, $p\ne2$ by Observation and  Lemma~\ref{deg=1} (a).  Suppose that $p=3$. By Lemma~\ref{deg=2}(a), none of  $x_1, x_2, x_3$ is a fair or bad vertex. Thus  all of $x_1, x_2, x_3$ must be good, contrary to Lemma~\ref{deg=1}(b).  Hence $p\ge4$. By Lemma~\ref{deg=2}(e), none of the vertices of $P$ are bad. If $x_2$ is fair,  then both  $x_1$ and $x_3$ must be good because neither $x_1$ nor $x_3$ are bad. By  Lemma~\ref{deg=1}(c) applied to $x_2$, $d_{G^*}(x_4)=3$, a contradiction. Thus  $x_2$ is good. Similarly,   $x_3$ is  good. By Lemma~\ref{deg=1}(g,\,c),  $x_1$ is neither good nor fair, and so  $x_1$ must be bad, a contradiction.   Thus  every component of $G^*[B]$ is isomorphic to $K_1$, $K_2$ or $P_3$. \medskip

 We shall apply the discharging method  to obtain a contradiction.
For each vertex $v\in V(G^*)$, let $\omega(v):= d_{G^*}(v)-\frac{24}{11}$ be the initial charge of $v$. Then $ \sum_{v\in V(G^*)} \omega(v) =2e(G^*)-\frac{24}{11}|G^*|=(2n_2+3n_3-n_1)-\frac{24}{11}(n_2+n_3)=\frac{9n_3-2n_2-11n_1}{11}<0$,  because $9n_3<2n_2+11n_1$.  Notice that for each $v\in V(G^*)$,  $\omega(v)=2-24/11=-2/11$ if  $d_{G^*}(v)=2$  and $\omega(v)=3-24/11=9/11$ if $d_{G^*}(v)=3$. We will redistribute the charges of vertices in $G^*$ according to the following discharging rule:\medskip

\noindent{\bf  (R)}: For each $x\in V(G^*)$ with $d_{G^*}(x)=3$ and  exactly $t\ge1$ neighbors of degree two in $G^*$,   $x$ sends $\frac{9}{11t} \ge \frac{3}{11}$ charges to each of its neighbors of degree two in $G^*$.\medskip

 Let $\omega^*$ be the new charge of  $G^*$ after applying the above discharging rule. We see that for any $v\in V(G^*)$ with $d_{G^*}(v)=3$, $\omega^*(v)\ge0$. We next show that  $\omega^*(B) :=\sum_{v\in B}\, \omega^*(v)\ge0$.
By the above claim,  each component $P$ of $G^*[B]$  is isomorphic to $K_1$, $K_2$ or $P_3$. Thus each endpoint of $P$ (with the endpoint of $P_1$ counted twice)  receives at least $3/11$ charge from its  neighbor in $V(G^*)\less B$ and so $\omega^*(P):=\sum_{v\in V(P)}\omega^*(v)\ge \sum_{v\in V(P)}\omega(v)+3/11+3/11\ge -6/11+3/11+3/11=0$.  Hence $\omega^*(B)=\sum_{P \in G^*[B]} \omega^*(P)\ge0$, where $P \in G^*[B]$ denotes that $P$ is a component of $G^*[B]$.  We see that
 $ \sum_{v\in V(G^*)} \omega^*(v) =\sum_{v\in V(G^*)\less B} \omega^*(v)+\omega^*(B)\ge0$, contrary to  the fact that $ \sum_{v\in V(G^*)} \omega^*(v)=\sum_{v\in V(G^*)} \omega(v) <0$.   \hfill\vrule height3pt width6pt depth2pt\bigskip

\noindent{\bf Remark.}   Kerdjoudj,  Kostochka and Raspaud~\cite{KKP2017} considered the list version of star edge-colorings of simple graphs. They proved that every subcubic graph is   star list-$8$-edge-colorable, and further proved the following.

\begin{thm} [\cite{KKP2017}]\label{KKP} Let $G$ be a subcubic simple graph.\medskip

\myitemitem {(a)} If $mad(G)<7/3$, then $G$ is  star list-$5$-edge-colorable.\medskip

\myitemitem {(b)} If $mad(G)<5/2$, then $G$ is  star list-$6$-edge-colorable.

\end{thm}

  \noindent {\bf Acknowledgments.} The authors would like to thank one anonymous referee for many helpful comments, in particular, for pointing out that  the Heawood graph is  star $5$-edge-colorable and $K_4$ with one subdivided edge has star chromatic index $6$, and bringing Reference [3] to our attention.  Yongtang Shi would like to thank Bojan Mohar for  his helpful discussion and for  mentioning the complexity problem on  star edge-coloring  during his visit to  Simon Fraser University.    The authors would like to thank Rong Luo for his helpful comments.

\end{document}